\documentclass[11pt,english]{amsart}
\usepackage[a4paper]{geometry}
\usepackage{amssymb}

\makeatletter

\newcommand{\noun}[1]{\textsc{#1}}

\usepackage{latexsym}
\theoremstyle{plain}
\newtheorem{theo}{Theorem}[section]

\newtheorem*{conjecture}{Conjecture}
\theoremstyle{remark}
\newtheorem{strem}[theo]{Remark}
\numberwithin{equation}{section} 
\numberwithin{figure}{section} 
\numberwithin{table}{section} 

\usepackage[bookmarksopen, naturalnames]{hyperref}

\makeatother

\usepackage{babel}

\begin{document}
\title{On a vector-valued Hopf-Dunford-Schwartz lemma}
\author{S. Charpentier, L. Deleaval}
\address{St\'ephane Charpentier, 
Laboratoire Paul Painlev\'e, UMR 8524, Universit\'e Lille 1, Cit\'e Scientifique, 59650 Villeneuve d'Ascq}
\email{Stephane.Charpentier@math.univ-lille1.fr}
\address{Luc Deleaval, Institut de Math\'ematiques de Jussieu, Universit\'e Pierre et Marie Curie, 4, Place Jussieu, 75005 Paris, France}
\email{deleaval@math.jussieu.fr}
\keywords{Hopf-Dunford-Schwartz lemma, semi-group, Dunkl analysis}
\subjclass[2010]{47A35, 42B25}

\begin{abstract}
In this paper, we state as a conjecture a vector-valued Hopf-Dunford-Schwartz lemma and give a partial answer to it. As an application of this powerful result, we prove some Fefferman-Stein inequalities in the setting of Dunkl analysis where the classical tools of real analysis cannot be applied. 
\end{abstract}

\maketitle

\section{Introduction}

A fundamental object in real and harmonic analysis is the Hardy-Littlewood maximal operator. Originally introduced by Hardy and Littlewood for functions defined on the circle (\cite{HardyLittlewood}), it was later extended for homogeneous spaces, and even in some non-homogeneous setting like noncompact symmetric spaces (\cite{stro}). Moreover, Fefferman and Stein extended it in a vector-valued setting and showed the following generalization of the Hardy-Littlewood theorem (\cite{FeffermanStein}), where we denote by $M$ the Hardy-Littlewood maximal operator given by 
\[
Mf(x)=\sup_{r>0} \frac{1}{|B_r(x)|}\int_{B_r(x)}|f(y)|\,\mathrm{d}y, \quad x \in \mathbb R^d,
\]
where $|X|$ is the Lebesgue measure of $X$ and $B_r(x)$ is the Euclidean ball of radius $x$ centered at the origin.

\begin{theo}\label{FSt}
Let  $f=(f_n)_{n \in \mathbb N}$ be a sequence of measurable functions on $\mathbb R^d$ and let $\overline{M}$ be the Fefferman-Stein operator given by $\overline{M}f=(Mf_n)_{n \in \mathbb N}$.
 Let $1<q<+\infty$.
\begin{enumerate}
\item If $\|f\|_{\ell^q} \in L^1(\mathbb R^d)$, then  for every $\lambda>0$ we have
\[
\Bigl|\Bigl\{x \in \mathbb R^d: \|\overline{M}f(x)\|_{\ell^q}>\lambda\Bigr\} \Bigr|\leqslant \frac{C}{\lambda}\bigl\|\|\,f\|_{\ell^q}\bigr\|_{L^1(\mathbb R^d)},
\]
where $C=C(q)$ is independent of $(f_n)_{n \in \mathbb N}$ and $\lambda$.
\item If  $\|f\|_{\ell^q} \in L^p(\mathbb R^d)$, with $1<p<+\infty$, then we have 
\[
\bigl\|\|\,\overline{M}f\|_{\ell^q}\bigr\|_{L^p(\mathbb R^d)}\leqslant C\bigl\|\|\,f\|_{\ell^q}\bigr\|_{L^p(\mathbb R^d)},
\]
where $C=C(p,q)$ is independent of $(f_n)_{n \in \mathbb N}$.
\end{enumerate}
\end{theo}

In contrast with the scalar case where different proofs may be given (combinatorics, covering method together with interpolation, operator theory with the scalar Hopf-Dunford-Schwartz lemma,...) there exists essentially one proof of this result (closely related to singular integral theory), based on three powerful tools of real analysis: the Hardy-Littlewood maximal theorem, a suitable Calder\'on-Zygmund decomposition and weighted inequalities. However, it seems to be important to give another proof, relying on a vector-valued version of the Hopf-Dunford-Schwartz lemma. Indeed, such a lemma should be very helpful in abstract setting, like Dunkl analysis, where a theory of singular integral is out of reach for the moment. More precisely, consider the Dunkl maximal operator $M_\kappa^W$ defined by
\begin{equation} \label{dunklmax}
M_\kappa^Wf(x)=\sup_{r>0}\frac{1}{\mu_\kappa^W(B_r)}\biggl|\int_{\mathbb R^d}f(y)\tau_x^{W}(\chi_{{}_{B_r}})(-y)\,\mathrm{d}\mu_\kappa^W(y)\biggr|, \quad x \in \mathbb R^d,
\end{equation}
where $\chi_{{}_{B_r}}$ is the characteristic function of the Euclidean ball of radius $r$ centered at the origin, $\tau_x^{W}$ is the Dunkl translation and $\mu_\kappa^W$ is a weighted Lebesgue measure invariant under the action of a reflection group $W$ (see Section 3 for more details). This operator, which reduces to the Hardy-Littlewood maximal operator in the case where the multiplicity function $\kappa$ is equal to $0$, is of particular interest in view of developing harmonic analysis associated with root systems and reflection groups. Nevertheless, the structure of the Dunkl translation prevents us from using the tools of real analysis (covering lemma, weighted inequalities...) and makes it difficult to study $M_\kappa^W$, even if the measure $\mu_\kappa^W$ is doubling. Mention that a scalar maximal theorem has been proved by use of the scalar Hopf-Dunford-Schwartz lemma (see \cite{Luc2,ThangaXuMaximal}), but no vector-valued extension for general $W$ is currently possible by using this result together with the other ingredients cited above. Nevertheless, Fefferman-Stein type inequalities have been established in \cite{Luc} when $W$ is isomorphic to $\mathbb Z_2^d$.  This was done thanks to explicit formulas in this particular case, which allow to bypass the lack of information on the Dunkl translation and so, to constuct a more convenient operator of Hardy-Littlewood type which controls pointwise $M_\kappa^{\mathbb Z_2^d}$.

The paper is organized as follows. In the next section, we make a brief review on the classical Hopf-Dunford-Schwartz lemma and we give a partial improvement of it in the vector-valued setting. Moreover, we state as a conjecture the general expected result. Section 3 is devoted to an application in Dunkl analysis.

Throughout the paper, we use $X\lesssim Y$ to denote the estimate $X\leqslant CY$ for some absolute constant $C$; if we need $C$ to depend on parameters, we indicate this by subscripts, thus $X\lesssim_p Y$ means that $X\leqslant C(p)Y$ for some $C(p)$ depending on $p$.

\section{A vector-valued Hopf-Dunford-Schwartz lemma}

\subsection{Background on Hopf-Dunford-Schwartz lemma}

Let $\left(\Omega,m\right)$ be a positive measure space. We
say that a one-parameter semi-group of linear operators $\left(T_{t}\right)_{t\geqslant0}$
satisfies the contraction property if each $T_{t}$ is a contraction
from $L^{1}\left(\Omega\right)$ into itself (for convenience, we
will just say $T_{t}$ is a contraction in $L^{1}\left(\Omega\right)$)
and in $L^{\infty}\left(\Omega\right)$. By Riesz convexity theorem,
$T_{t}$ is also a contraction in $L^{p}\left(\Omega\right)$, for
any $1\leqslant p\leqslant+\infty$ and any $t\geqslant0$. Moreover, if such a semi-group
$\left(T_{t}\right)_{t\geqslant0}$ is assumed to be strongly measurable, then
$\left(T_{t}\right)_{t\geqslant0}$ is strongly integrable over every finite
interval, i.e. for any $f\in L^{p}\left(\Omega\right)$, $t\mapsto T_{t}\left(f\right)$
is integrable with respect to the Lebesgue measure along every interval
$0\leqslant t\leqslant\alpha$. This allows to consider the averages\[
A_{\alpha}f:=\frac{1}{\alpha}\int_{0}^{\alpha}T_{t}\left(f\right)\,\mathrm{d}t,\quad \alpha>0.\]
We set $A_{0}\left(f\right)=f$ for any $f\in L^{p}\left(\Omega\right)$.
Because the linear map $f\mapsto A_{\alpha}f$ is closed, $A_{\alpha}$
is a contraction in each $L^{p}\left(\Omega\right)$, $1\leqslant p\leqslant+\infty$.
In particular, the function $x\mapsto\frac{1}{\alpha}\int_{0}^{\alpha}T_{t}\left(f\right)\left(x\right)\,\mathrm{d}t$
is $m$-measurable as well as the maximal function $x\mapsto M_Tf\left(x\right)$
defined by\begin{equation}
M_Tf\left(x\right):=\sup_{\alpha\geqslant0}\frac{1}{\alpha}\left|\int_{0}^{\alpha}T_{t}\left(f\right)(x)\,\mathrm{d}t\right|,\label{eq|def-max-funct-scalar}\end{equation}
for every $x\in\Omega$ except those in a $m$-null subset of $\Omega$.
For all the details, we refer to \cite[Chapter VIII.7]{Dunford}.

The scalar Hopf-Dunford-Schwartz ergodic lemma may be stated as follows:
\begin{theo}
[Lemma VIII.7.6 and Theorem VIII.7.7 of \cite{Dunford}]\label{thm|scalar-Hopf-lemma}Let
$\left(T_{t}\right)_{t\geqslant0}$ be a strongly measurable semi-group which
satisfies the contraction property. Let $f$ be a measurable function on $\Omega$. \begin{enumerate}\item If $f \in L^1(\Omega)$, then for
every $\lambda>0$ we have\[
{\displaystyle m\left(\left\{ x\in\Omega:  M_Tf\left(x\right)>\lambda\right\} \right)\leqslant\frac{2}{\lambda}\left\Vert f\right\Vert _{L^{1}(\Omega)}.}\]
\item If $f \in L^p(\Omega)$, with $1<p<+\infty$, then we have\[
\left\Vert M_Tf\right\Vert _{L^{p}(\Omega)}\leqslant2\left(\frac{p}{p-1}\right)^{1/p}\left\Vert f\right\Vert _{L^{p}(\Omega)}.\]
\end{enumerate}
\end{theo}

The first version of this theorem goes back to E. Hopf \cite{EHopf} and concerned ergodic means. A generalization was given by Dunford and Schwartz in \cite{Dun-Schwartz,Schwartz} leading to the above statement. We also refer to \cite{Krengel} for a nice presentation
of ergodic theorems.

As mentioned in \cite{Stein}, the Hopf-Dunford-Schwartz lemma is a powerful tool in harmonic
analysis. For instance, it allows to prove in abstract setting the boundedness of some maximal operators as soon as the relationship between averages on balls and the heat semi-group (or Poisson semi-group) could be exploited. Because difficulties arise when formulating estimates for operators in the context of vector-valued functions, an extension of an abstract tool such as Hopf-Dunford-Schwartz lemma would be very useful.

Actually, a lot of vector-valued versions of this result have been given by
various authors (for e.g. \cite{Chacon,Taggart,Yoshimoto} and the references
therein). All of them consist in generalizing the maximal function
operator (\ref{eq|def-max-funct-scalar}) by replacing the absolute
value by a norm in a Banach space. Different approaches to this
problem were considered: for instance, in \cite{Chacon}, the author gave a tricky
but specific proof to obtain such a vector-valued Hopf-Dunford-Schwartz lemma (his
result is not stated in terms of averages of semi-group but in terms of ergodic means,
yet one can obtain it in full generality arguing as in \cite{Dunford})
without requiring any results allowing him to consider positive
operator. Recently, Taggart \cite{Taggart} gave a direct proof for
$1<p<+\infty$ appealing to classical Banach-valued extension results
for positive operators (see \cite{Rubio,Grafakos}). Yet,
this positivity assumption requires some non-trivial background.

Let us state, in a slightly different form, the vector-valued Hopf-Dunford-Schwartz lemma obtained in \cite{Chacon,Taggart}:
\begin{theo}
Assume that $B$ is a Banach space and that $\left(T_{t}\right)_{t\geqslant0}$
is a strongly continuous semi-group which satisfies the contraction property. Define\[
\widetilde{M_T}f:=\sup_{\alpha\geqslant0}\frac{1}{\alpha}\left\Vert \int_{0}^{\alpha}\widetilde{T_{t}}\left(f\right)\,\mathrm{d}t\right\Vert _{B},\]
where $\widetilde{T_{t}}$ is the linear extension of $T_{t}$ to
the Banach space $L^{p}\left(\Omega,B\right)$ of vector-valued functions
$f$ such that $x\mapsto\left\Vert f(x)\right\Vert _{B}$
is measurable on $\Omega$ and belongs to $L^{p}\left(\Omega\right)$.
 \begin{enumerate}\item If $f\in L^{1}\left(\Omega,B\right)$, then for every $\lambda>0$ we have\[
m\left(\left\{ x\in\Omega: \,\widetilde{M_T}f\left(x\right)>\lambda\right\} \right)\leqslant\frac{2}{\lambda}\left\Vert f\right\Vert _{L^{1}\left(\Omega,B\right)}.\]
\item If $f\in L^{p}\left(\Omega,B\right)$, with $1<p<+\infty$, then we have \[
\left\Vert \widetilde{M_T}f\right\Vert _{L^{p}\left(\Omega,B\right)}\leqslant2\left(\frac{p}{p-1}\right)^{1/p}\left\Vert f\right\Vert _{L^{p}\left(\Omega,B\right)}.\]
\end{enumerate}
\end{theo}

Because of the definition of the maximal operator $\widetilde{M_T}$ under consideration, this result is useless to obtain, in particular, Fefferman-Stein type inequalities. In the next paragraph, we introduce another extension which will fill this lack.

\subsection{A partial improvement of Hopf-Dunford-Schwartz lemma}

As in the previous paragraph, $M_T$ stands for the maximal operator
defined by (\ref{eq|def-max-funct-scalar}). We consider the Banach
space $L^{p}\left(\Omega,l^{q}\right)$ ($L^{p}\left(l^{q}\right)$ for short),
$1\leqslant p,q\leqslant+\infty$, consisting of all $l^{q}$-valued functions
$f=\left(f_{n}\right)_{n\in\mathbb{N}}$ (each $f_{n}$ is a measurable
real or complex valued function) for which $x\mapsto\left(\sum_{n=0}^{+\infty}\left|f_{n}\left(x\right)\right|^{q}\right)^{1/q}$
is finite $m$-a.e. and such that\[
\left\Vert f\right\Vert _{L^{p}\left(l^{q}\right)}:=\left(\int_{\Omega}\left(\sum_{n=0}^{+\infty}\left|f_n\left(x\right)\right|^{q}\right)^{p/q}\,\mathrm{d}m\left(x\right)\right)^{1/q}<+\infty\]
(with the usual modification for $p$ or $q=+\infty$). Now, for $f=\left(f_{n}\right)_{n \in \mathbb N}\in L^{p}\left(l^{q}\right)$,
we define the following vector-valued maximal function operator
$$\overline{M_T}f:=\left(M_Tf_{n}\right)_{n\in\mathbb{N}}=\left(\sup_{\alpha\geqslant0}\frac{1}{\alpha}\left|\int_{0}^{\alpha}T_{t}\left(f_{n}\right)\,\mathrm{d}t\right|\right)_{n\in\mathbb{N}}.$$

It is now very natural to wonder whether the conclusions of the original
scalar Hopf-Dunford-Schwartz lemma are also true for $\overline{M_T}$ when replacing
$L^{p}\left(\Omega\right)$ by $L^{p}\left(l^{q}\right)$, $1\leqslant p,q\leqslant+\infty$.
The main theorem of this paragraph is a partial positive answer to
that question, namely for $1<p\leqslant q<+\infty$:
\begin{theo}
\label{thm|Hopf-vector-valued}Let $\left(T_{t}\right)_{t\geqslant0}$ be a
strongly continuous semi-group which satisfies the contraction property. Assume that $1<p\leqslant q<+\infty$.
Then, for any $f\in L^{p}\left(l^{q}\right)$ we have\begin{equation}
\left\Vert \overline{M_T}f\right\Vert _{L^{p}\left(l^{q}\right)}\lesssim_{p,q}\left\Vert f\right\Vert _{L^{p}\left(l^{q}\right)}.\label{eq|Hopf-ineq-vector}\end{equation}

\end{theo}

\begin{strem}
The case $p=q$ is a straightforward consequence of the scalar case. Because Theorem \ref{thm|scalar-Hopf-lemma} holds for $p=+\infty$, one can easily check that Inequality (\ref{eq|Hopf-ineq-vector}) also holds for $1<p<+\infty$ and $q=+\infty$.
\end{strem}

The first step consists in showing that we can reduce the problem
to a semi-group of linear \emph{positive} operators $S_{t}$. Recall that
an operator $T$ on $L^{p}\left(\Omega\right)$ (non necessarily linear) is positive if $\left|T\left(f\right)\right|\leqslant T\left(g\right)$
whenever $\left|f\right|\leqslant g$ almost everywhere. This is a direct combination of the proof of \cite[Theorem 3.1]{Taggart}
and \cite[Lemma 3.5]{Taggart}:
\begin{theo}
Suppose that $\left(T_{t}\right)_{t\geqslant0}$ is a strongly continuous semi-group
which satisfies the contraction
property. Then there exists a strongly continuous semi-group $\left(S_{t}\right)_{t\geqslant0}$
satisfying the contraction property, such that $S_{t}$ is linear positive for any $t\geqslant0$, and which dominates $T$ in the sense that $$\left| T_{t}f\right| \leqslant S_{t}|f|,\quad \forall f\in L^{p}\left(\Omega\right),\, \forall t\geqslant0.$$
\end{theo}
Therefore, to show Theorem \ref{thm|Hopf-vector-valued}, we are reduced to prove that (\ref{eq|Hopf-ineq-vector}) is satisfied by the
maximal operator $\overline{M_T^{+}}$ instead of $\overline{M_{T}}$, for any $f=\left(f_{n}\right)_{n\in \mathbb N}\in L^{p}\left(l^{q}\right)$ positive (i.e. $f_{n}\geqslant0$, $\forall n\in \mathbb{N}$) where\[
\overline{M_T^{+}}f:=\left(\sup_{\alpha\geqslant0}\frac{1}{\alpha}\int_{0}^{\alpha}S_{t}\left(f_{n}\right)\,\mathrm{d}t\right)_{n\in\mathbb{N}},\]
with $\left(S_{t}\right)_{t\geqslant0}$ given by the previous theorem. Moreover,
if we denote by $M_S^{+}$ the scalar maximal operator $M_S^{+}f=\sup_{\alpha\geqslant0}\frac{1}{\alpha}\int_{0}^{\alpha}S_{t}\left(f\right)\,\mathrm{d}t$,
$f\in L^{p}\left(\Omega\right)$, then $M_S^{+}$ is linearizable in
the sense of \cite[Definition V.1.20]{Rubio}, that is\[
M_S^{+}=\left\Vert A_{+}\right\Vert _{L^{\infty}\left(\mathbb{R}\right)},\]
where $A_{+}$ is the linear operator in $L^{p}\left(\Omega\right)$
which takes $f$ to $\left(\alpha\mapsto{\displaystyle \frac{1}{\alpha}}\int_{0}^{\alpha}S_{t}\left(f\right)\,\mathrm{d}t\right)$.
In addition, by 2) of the scalar Hopf-Dunford-Schwartz lemma (Theorem \ref{thm|scalar-Hopf-lemma}),
$M_S^{+}$ is bounded in $L^{p}\left(\Omega\right)$ for any $p>1$.

The proof of Theorem \ref{thm|Hopf-vector-valued} is then a consequence
of the following:
\begin{theo}[Corollary V.1.23 of \cite{Rubio}]\label{thm|Rubio-linearizable}Let $T$ be a linearizable operator
which is bounded in $L^{p}\left(\Omega\right)$ for some $1\leqslant p<+\infty$.
If $T$ is positive, then for any $q\geqslant p$, $T$ has a bounded extension
to $L^{p}\left(l^{q}\right)$, i.e.\[
\left\Vert \left(\sum_{n=0}^{+\infty}\left|T\left(f_{n}\right)\right|^{q}\right)^{1/q}\right\Vert _{L^{p}\left(\Omega\right)}\lesssim_{p,q}\left\Vert \left(\sum_{n=0}^{+\infty}\left|f_{n}\right|^{q}\right)^{1/q}\right\Vert _{L^{p}\left(\Omega\right)},\]
for any $f=\left(f_{n}\right)_{n \in \mathbb N}$ in $L^{p}\left(l^{q}\right)$.
\end{theo}
\medskip{}

\begin{strem} (1) Theorem \ref{thm|Hopf-vector-valued} allows to recover the Fefferman-Stein
inequalities of Theorem \ref{FSt} when $1<p\leqslant q$. Indeed, it is a simple consequence of the pointwise inequality \[Mf(x)\lesssim M_{H}f(x),\quad x\in\mathbb R^{d},\] where $M_{H}$ is the standard Euclidean heat maximal operator.

\noindent{} (2) Let us observe that Theorem \ref{thm|Hopf-vector-valued} can still be extended for the maximal operator
\begin{equation}
\overline{M_{\widetilde{T}}^{B}}f:=\left(M_{\widetilde{T}}^{B}f_{n}\right)_{n\in \mathbb{N}}=\left(\sup_{\alpha\geqslant0}\frac{1}{\alpha}\left\Vert \int_{0}^{\alpha}\widetilde{T_{t}}\left(f_{n}\right)\,\mathrm{d}t\right\Vert _{B}\right)_{n\in\mathbb{N}},\end{equation}
where $B$ is a Banach space, $\widetilde{T_{t}}$ is the linear extension of $T_{t}$ to $L^{p}\left(\Omega,B\right)$ and $\left(f_{n}\right)_{n\in \mathbb{N}}$ is a sequence of $B$-valued functions. Indeed, it suffices to follow the proof of Theorem \ref{thm|Hopf-vector-valued}, by using Taggart's result \cite[Corollary 4.2]{Taggart} and by noticing that, if $\left(S_{t}\right)_{t\geqslant0}$ is a positive semi-group dominating $\left(T_{t}\right)_{t\geqslant0}$, then $M_{\widetilde{T}}^{B,+}:=M_{\widetilde{S}}^{B}$ is a linearizable positive operator, that is for $f\in L^{p}\left(\Omega,B\right)$,
\begin{equation}{M_{\widetilde{T}}^{B,+}}f=\left\Vert \left(\alpha\mapsto{\displaystyle \frac{1}{\alpha}}\int_{0}^{\alpha}\widetilde{S}_{t}\left(f\right)\,\mathrm{d}t\right)\right\Vert _{L^{\infty}\left(\mathbb{R},B\right)},\end{equation}
which is bounded in $L^{p}\left(\Omega,B\right)$.
\noindent{}In particular, such a further extension of Theorem \ref{thm|Hopf-vector-valued} is a (self-) improvement of Taggart's result \cite[Corollary 4.2]{Taggart}.
\end{strem}

Since Theorem \ref{thm|Rubio-linearizable} does not hold for arbitrary $1\leqslant p,q\leqslant +\infty$ (see \cite{Rubio}), our method does not give a complete vector-valued version of Hopf-Dunford-Schwartz lemma. Yet, we make the following conjecture, a confirmation of which would be a powerful result.
\begin{conjecture}
Let $\left(T_{t}\right)_{t\geqslant0}$ be a strongly continuous semi-group which satisfies the
contraction property. Let $1<q<+\infty$. \begin{enumerate}\item If $f\in L^{1}\left(l^{q}\right)$, then for every $\lambda>0$ we have
\[
{\displaystyle m\left(\left\{ x\in\Omega: \left\Vert \overline{M_T}f\left(x\right)\right\Vert _{l^{q}}>\lambda\right\} \right)\lesssim_{q}\frac{1}{\lambda}\left\Vert f\right\Vert _{L^{1}\left(l^{q}\right)}.}\]
\item If $f\in L^{p}\left(l^{q}\right)$, with $1<p<+\infty$, then we have \[
\left\Vert \overline{M_T}f\right\Vert _{L^{p}\left(l^{q}\right)}\lesssim_{p,q}\left\Vert f\right\Vert _{L^{p}\left(l^{q}\right)}.\]
\end{enumerate}
\end{conjecture}

\section{Some Fefferman-Stein inequalities for the Dunkl maximal operator}
In this section, we present an application of our vector-valued Hopf-Dunford-Schwartz lemma (Theorem \ref{thm|Hopf-vector-valued}) in the setting of Dunkl analysis, which extends Fourier analysis on Euclidean spaces and analysis on Riemannian symmetric spaces of Euclidean type. Before stating our Fefferman-Stein type inequalities for the operator $\overline{M^W_\kappa}f=\left(M^W_\kappa f_n\right)_{n \in \mathbb N}$ (where $M_\kappa^W$ is given by \eqref{dunklmax}), we give a brief account on the Dunkl theory for the reader's convenience.

\subsection{Background on Dunkl analysis}
For a large panorama of this theory, we refer to \cite{Roslerlecturenotes} and the references therein.

Let $W \subset \mathcal O(\mathbb{R}^d)$ be a finite reflection group associated with a reduced root system $\mathcal R$ (not necessarily crystallographic) and let $\kappa:\mathcal R\to \mathbb{R}_+$ be a multiplicity function, that is a $W$-invariant function.

The (rational) Dunkl operators $T_\xi^{\mathcal R}$ on $\mathbb R^d$, which were introduced in \cite{Dunkl89}, are the following $\kappa$-deformations of directional derivatives $\partial_\xi$ by reflections
\[
T_\xi^{\mathcal R}f(x)=\partial_\xi f(x)+\frac{1}{2}\sum_{\alpha \in \mathcal R}\kappa(\alpha)\frac{f(x)-f( \sigma_\alpha(x))}{\langle x,\alpha\rangle}\langle\xi,\alpha\rangle,\quad x \in \mathbb R^d,
\]
where $\sigma_\alpha$ denotes the reflection with respect to the hyperplane orthogonal to $\alpha$. The most important property of these operators is their commutativity  (\cite{Dunkl89}). Therefore, we are naturally led to consider the eigenfunction problem 
\begin{equation} \label{spectral} T^{\mathcal R}_\xi f=\langle y,\xi\rangle f, \quad \forall \xi \in \mathbb R^d,\end{equation}
with $y \in \mathbb C^d$ a fixed parameter. Opdam has completely solved this problem (\cite{Opdam}).

\begin{theo}
Let $y \in \mathbb C^d$. There exists a unique $f=E_\kappa^W(\cdot,y)$ solution of \eqref{spectral} which is real-analytic on $\mathbb R^d$ and satisfies $f(0)=1$. Moreover $E_\kappa^W$, called the Dunkl kernel, extends to a holomorphic function on $\mathbb C^d\times \mathbb C^d$.
\end{theo} 
Unfortunately, the Dunkl kernel is explicitly known only in very few cases; nevertheless we know that this kernel has many properties in common with the classical exponential to which  it reduces when $\kappa=0$. The Dunkl kernel is of particular interest as it gives rise to an integral transform, which is taken with respect to a weighted Lebesgue measure invariant under the action of $W$ and which generalizes the Euclidean Fourier transform.

More precisely, we introduce the homogeneous (of degree $\gamma$) and $W$-invariant measure \[\mathrm{d}\mu_\kappa^W(x)=\prod_{\alpha \in \mathcal R}|\langle x,\alpha\rangle|^{\kappa(\alpha)}\,\mathrm{d}x.\] We denote  by $L^p_\kappa$ the space $L^p(\mathbb R^d; \mu_\kappa^W)$ (for $1\leqslant p\leqslant +\infty$) and we use the shorter notation $\mathopen\|\cdot\mathclose\|_{\kappa,p}$ instead of $\mathopen\|\cdot\|_{L_\kappa^p}$.
Then for every $f \in L^1_{\kappa}$, the Dunkl transform of $f$, denoted by $\mathcal F_\kappa^W(f)$, is defined by
\[
\mathcal F_\kappa^W(f)(x)=c_\kappa^W\int_{\mathbb R^d}E^W_\kappa(-ix,y)f(y)\,\mathrm{d}\mu_\kappa^W(y), \quad x \in \mathbb R^d, 
\]
where $c_\kappa^W$ is a Mehta-type constant. Let us point out that the Dunkl transform coincides with the Euclidean Fourier transform when $\kappa=0$ and that it is more or less a Hankel transform when $d=1$. The two main properties of the Dunkl transform are given in the following  theorem (\cite{deJeuInv,DunklHankel}).
\begin{theo}\emph{(1)} \textbf{Inversion formula.} Let $f \in L^1_\kappa$. If $\mathcal F_\kappa^W(f)$ is in $L^1_\kappa$, then we have the following inversion formula
\[
 f(x)=c_\kappa^W\int_{\mathbb R^d}E^W_\kappa(ix,y)\mathcal F_\kappa^W(f)(y)\,\mathrm{d}\mu_\kappa^W(y).
\]
\noindent{}\emph{(2)} \textbf{Plancherel theorem.} The Dunkl transform has a unique extension to an isometric isomorphism of $L^2_\kappa$.
\end{theo}

The Dunkl transform shares many other properties with the Fourier transform. Therefore, it is natural to associate a generalized translation operator with this transform.

There are many ways to define the Dunkl translation but we use the definition which most underlines the analogy with the Fourier transform. It is the definition given in \cite{ThangaXuMaximal} with a different convention. Let $x \in \mathbb R^d$. The Dunkl translation  $f\mapsto\tau_x^{W}f$ is defined on $L^2_{\kappa}$ by the equation
\[
\mathcal F_\kappa^W(\tau_x^{W}f)(y)=E_\kappa^W(ix,y)\mathcal F_\kappa^W(f)(y),\quad y \in \mathbb R^d.
\]
In Fourier analysis, the translation operator $f\mapsto f(\cdot+x)$ (to which the Dunkl translation reduces when $\kappa=0$) is positive and $L^p$-bounded. In the Dunkl setting, $\tau_x^W$ is not a positive operator (\cite{RoslerHyp,ThangaXuMaximal}) and the $L^p_\kappa$-boundedness is still a challenging problem, apart from the trivial case where $p=2$ (thanks to the Plancherel theorem and the fact that $|E_\kappa^W(ix,y)|\leqslant1$). Moreover, its structure prevents from using covering methods.

\subsection{Fefferman-Stein type inequalities for the Dunkl maximal operator}

With all these definitions in mind, we can now state our main result.

\begin{theo}
\label{thm|Hopf-Dunkl}Let $1<p\leqslant q<+\infty$.
Then, for any $f\in L_\kappa^{p}\left(l^{q}\right)$ we have\begin{equation}
\left\Vert \overline{M_\kappa^W}f\right\Vert _{L^p_\kappa\left(l^{q}\right)}\lesssim_{p,q}\left\Vert f\right\Vert _{L^{p}_\kappa\left(l^{q}\right)}.\label{eq|Hopf-ineq-vector-final}\end{equation}
\end{theo}

\begin{strem} In \cite{Luc}, complete Fefferman-Stein inequalites were given in the particular case where $W\simeq\mathbb Z_2^d$ but, as already claimed, the proof relies on some explicit formulas for $\tau_x^{\mathbb Z_2^d}$ which allow to construct a $\mathbb Z_2^d$-invariant maximal operator, which controls pointwise $M_\kappa^{\mathbb Z_2^d}$ and for which we can use covering argument together with interpolation. 
\end{strem}

We come to the proof of Theorem \ref{thm|Hopf-Dunkl}. Let us introduce the so-called Dunkl-type heat semi-group $(H_t^{W})_{t\geqslant0}$ which is associated with the Dunkl Laplacian $\Delta_\kappa^W=\sum_{j=1}^d(T_{e_j}^{\mathcal R})^2$ (see \cite{RoslerHermite}).  More precisely for every $f \in L^p_\kappa$, with $1\leqslant p \leqslant+\infty$, and for every $t\geqslant0$, it is given by
\[H_t^{W}f=\begin{cases}\int_{\mathbb R^d}f(y)Q_\kappa^W(\cdot,y,t)\,\mathrm{d}\mu_\kappa^W(y) &\text{if}\ t>0
\\ f &\text{if}\ t=0,
 \end{cases}\]
 where \[
Q_\kappa^W(x,y,t)=\frac{c_\kappa^W}{(2t)^{\frac{d}{2}+\gamma}}\mathrm{e}^{-\frac{(\|x\|^2+\|y\|^2)}{4t}}E_\kappa^W\Bigl(\frac{x}{\sqrt{2t}},\frac{y}{\sqrt{2t}}\Bigr)>0,\quad x,y \in \mathbb R^d, t>0. 
\]
According to \cite{Luc2} (see also \cite{ThangaXuMaximal}), we have the following pointwise inequality
\[M_\kappa^Wf(x)\lesssim\sup_{\alpha>0}\frac{1}{\alpha}\int_0^\alpha H_t^{W}|f(x)|\,\mathrm{d}t:=M_{\kappa,H}^Wf(x), \quad x \in \mathbb R^d.\] Thanks to this inequality, it suffices to show that Inequality (\ref{eq|Hopf-ineq-vector-final}) is true by replacing $\overline{M^W_\kappa}f=\left(M^W_\kappa f_n\right)_{n \in \mathbb N}$ by $\overline{M_{\kappa,H}^W}f=\left(M_{\kappa,H}^W f_n\right)_{n \in \mathbb N}$, for any $f=\left(f_n\right)_{n \in \mathbb N}$ in $L_\kappa^{p}\left(l^{q}\right)$ positive. Now $(H_t)_{t\geqslant0}$ is a symmetric diffusion semi-group on $L_\kappa^p$ for $1\leqslant p \leqslant +\infty$ (Theorem  2.6. in \cite{Luc2}). Thus, by our vector-valued Hopf-Dunford-Schwartz lemma (Theorem \ref{thm|Hopf-vector-valued}), we get the desired estimate (\ref{eq|Hopf-ineq-vector-final}).

\end{document}